\documentclass[11pt]{amsart}

\usepackage[T1]{fontenc}
\usepackage{lmodern}
\usepackage{amsmath,amssymb,amsthm,mathtools}
\usepackage{microtype}
\usepackage[hidelinks]{hyperref}

\newtheorem{theorem}{Theorem}[section]
\newtheorem{proposition}[theorem]{Proposition}
\newtheorem{lemma}[theorem]{Lemma}

\newtheorem{remark}[theorem]{Remark}

\newcommand{\Area}{\operatorname{Area}}
\newcommand{\KS}{K_S}
\newcommand{\Kpos}{\mathcal K_{+}}
\newcommand{\XS}{X_S}
\newcommand{\UC}{U_{\mathrm C}}
\newcommand{\R}{\mathbb R}

\newcommand{\eps}{\varepsilon}
\newcommand{\ol}{\overline}

\title[prescribed combinatorial curvature beyond positivity]{Prescribed Combinatorial Curvature Beyond Positivity for Spherical Vertex Scalings on The Sphere}
\author{Chen Liu, Xu Xu}
\date{\today}

\address{School of Mathematics and Statistics, Wuhan University, Wuhan, 430072, P.R. China} \email{2023302021080@whu.edu.cn}

\address{School of Mathematics and Statistics, Wuhan University, Wuhan, 430072, P.R. China} \email{xuxu2@whu.edu.cn}



\begin{document}

\begin{abstract}
Let $(S,V)$ be a marked 2-sphere with $|V|\ge3$, and fix a spherical discrete conformal class.  
Izmestiev-Prosanov-Wu proved existence and uniqueness for positive prescribed combinatorial curvature in the spherical discrete conformal class on $(S,V)$.  
In this paper, we identify a region of prescribed combinatorial curvature vectors, allowing some negative components, for which existence still holds.
\end{abstract}

\maketitle

\section{Introduction}

Vertex scaling is a natural discretization of smooth conformal transformations in Riemannian geometry.
It was first introduced by Luo \cite{Luo1} in the Euclidean background geometry, and then by Bobenko–Pinkall–Springborn \cite{BPS15} in the hyperbolic and spherical background geometries.
The prescribed combinatorial curvature problem on surfaces has been one of the fundamental problems in the theory of vertex scalings.
In the Euclidean background geometry, Gu–Luo–Sun–Wu \cite{GLSW18} proved that the prescribed combinatorial curvature problem is solvable if and only if the prescribed combinatorial curvature satisfies the discrete Gauss–Bonnet equality.
In the hyperbolic background geometry, Gu–Guo–Luo–Sun–Wu \cite{GGLSW18} proved that the prescribed combinatorial curvature problem is solvable if and only if the prescribed combinatorial curvature satisfies the corresponding discrete Gauss–Bonnet inequality.
However, in the spherical background geometry, the prescribed combinatorial curvature problem is much more subtle.
Recently, Izmestiev-Prosanov-Wu \cite{IzmestievProsanovWu24} solved the prescribed combinatorial curvature problem for spherical vertex scaling on the sphere when the prescribed combinatorial curvature is positive. Izmestiev-Prosanov-Wu's result is a discrete analogue of the well-known result of Luo-Tian \cite{LuoTian92} and Troyanov \cite{Troyanov91}. In this paper, we identify a region of prescribed combinatorial curvature vectors, allowing some negative components, for which existence still holds.

Let $S$ be the two-sphere and let $V\subset S$ be a finite set.  A spherical cone metric on $(S,V)$ is locally isometric to the unit sphere away from $V$.  If $\omega_v$ denotes the total angle at $v$, its combinatorial curvature (sometimes discrete curvature or cone curvature) is
\[
  \kappa_v=2\pi-\omega_v<2\pi.
\]
The discrete Gauss-Bonnet formula for a spherical cone metric $d$ on the sphere is as follows:
\begin{equation}\label{eq:GB}
  \Area(d)+\sum_{v\in V}\kappa_v(d)=4\pi,
\end{equation}
where $Area(d)$ is the area of the sphere in the spherical cone metric $d$.

For a fixed smooth conformal structure on the sphere, Troyanov \cite{Troyanov91} and Luo-Tian \cite{LuoTian92} proved an existence and uniqueness theorem for spherical cone metrics with positive cone curvatures; the conditions are the discrete Gauss-Bonnet inequality 
$$\sum_{v\in V}\kappa_v(d)<4\pi$$
together with the inequality 
$$ 2\kappa_v<\sum_{w\in V}\kappa_w, \ \ \forall v\in V.$$  

In the spherical discrete setting considered here, Izmestiev-Prosanov-Wu \cite{IzmestievProsanovWu24} proved the corresponding discrete version of Troyanov and Luo-Tian's result in the same positive curvature region.  
In this paper, 
we establish realization on the region $K_S$, which allows some negative components.

\begin{theorem}\label{thm:main}
Let $(S,V)$ be a marked two-sphere with $n:=|V|\ge3$, and let $d$
be a spherical cone metric on $(S, V)$ admitting a geodesic triangulation by convex spherical triangles.  Suppose that
$\ol\kappa\in(-\infty,2\pi)^V$
satisfies
\begin{equation}\label{eq:main-conditions}
  2\pi<\sum_{v\in V}\ol\kappa_v<4\pi,
  \qquad
  2\ol\kappa_v<\sum_{w\in V}\ol\kappa_w
  \quad (v\in V).
\end{equation}
Then there exists a spherical cone metric $d'$ discretely conformal to $d$ with
$\kappa(d')=\ol\kappa$.
Furthermore, if $\ol\kappa$ is a regular value of the curvature map, then the number of such metrics is finite and odd.
\end{theorem}

We can compare Theorem \ref{thm:main} with the result of Izmestiev-Prosanov-Wu in \cite{IzmestievProsanovWu24}. 
Theorem~\ref{thm:main} keeps the same discrete Gauss-Bonnet inequality and inequality $2\kappa_v<\sum_{w\in V}\kappa_w$, while adding the lower bound $\sum_v\kappa_v>2\pi$ and removing the condition $\kappa_v>0$.  Thus the region in Theorem~\ref{thm:main} and the region in Izmestiev-Prosanov-Wu's result are not nested.  The new feature is that prescribed combinatorial curvatures with some components negative are allowed. For example, if $n\geq 4$, the prescribed combinatorial curvature n-vector
\[
  \left(-\eps,\frac{3\pi+\eps}{n-1},\frac{3\pi+\eps}{n-1}, \cdots, \frac{3\pi+\eps}{n-1}\right)
\]
satisfies \eqref{eq:main-conditions} whenever $0<\eps<\frac{3n-9}{2}\pi$.  
\begin{remark}
    Negative curvature components cannot occur when $|V|=3$ in Theorem \ref{thm:main}. Indeed, suppose that $\ol\kappa_1<0$, then by the inequality $2\ol\kappa_v<\sum_{w\in V}\ol\kappa_w$, we have 
    $$\ol\kappa_2<\ol\kappa_1+\ol\kappa_3,\qquad \ol\kappa_3<\ol\kappa_1+\ol\kappa_2.$$
    This would imply $\ol\kappa_2<\ol\kappa_3<\ol\kappa_2$, a contradiction.
\end{remark}
The proof of Theorem~\ref{thm:main} is based on an application of degree theory and a refinement of the arguments of Izmestiev-Prosanov-Wu in \cite{IzmestievProsanovWu24}.  The key point is a boundary estimate: a finite Delaunay degeneration has limiting area at least $2\pi$.  By \eqref{eq:GB}, the lower bound in \eqref{eq:main-conditions} excludes this degeneration.  Degree theory then gives existence without requiring global injectivity.  

Let \(U'(d)\subset\mathbb R^V\) be the parameter domain of Izmestiev–Prosanov–Wu, and let \(U(d)\) be the connected component of \(U'(d)\) containing the origin. Their parametrization identifies \(U(d)\) with \(\mathcal C_+(d)\).
The combinatorial curvature is a $C^1$ map
\begin{equation}\label{eq:kappa-map}
  \kappa:U(d)\longrightarrow(-\infty,2\pi)^V.
\end{equation}
Set
\begin{equation}\label{eq:KS}
  \KS:=\left\{\kappa\in\R^V:
  2\pi<\sum_{v\in V}\kappa_v<4\pi,
  \quad
  \quad 2\kappa_v<\sum_{w\in V}\kappa_w, \forall v\in V
  \right\},
\end{equation}
and
\begin{equation}\label{eq:XS}
  \XS(d):=\kappa^{-1}(\KS).
\end{equation}
We shall prove that the map $\kappa:\XS(d)\longrightarrow\KS$
is proper and has degree $\pm1$.  The oddness assertion in Theorem~\ref{thm:main} follows at once for regular values.

The paper is arranged as follows.
Section~2 recalls the notation and some results from \cite{IzmestievProsanovWu24}.  Section~3 contains the compactness lemmas.  Theorem~\ref{thm:main} is proved in Section~4.
\\
\\
\textbf{Acknowledgments} 
The research of Chen Liu is supported by the National Natural Science Foundation of China under Grant No. 125B1006.
The research of Xu Xu is supported by the National Natural Science Foundation of China under Grant No. 12471057 and Hubei Provincial Natural Science Fund for Outstanding Research
Groups under Grant No. 2026AFA008. 
\\
\\
\noindent\textbf{AI Statement:} Generative AI tools were used during the preparation of this manuscript to assist with English-language editing, figure editing and LaTeX typesetting. All mathematical statements and proofs were independently given by the authors, who take full responsibility for the content.

\section{Settings and known results}\label{sec:setting}

We use the notation from Izmestiev-Prosanov-Wu \cite{IzmestievProsanovWu24}, and recall some results in \cite{IzmestievProsanovWu24} that are needed in this paper.

\subsection{Cone metrics and Delaunay triangulations}

For a triangulation $T$ of $(S,V)$, let $\mathfrak M_+(T)$ be the set of spherical cone metrics realizing $T$ by convex spherical triangles.  In edge-length coordinates,
\[
  \mathfrak M_+(T)=\left\{\ell\in(0,\pi)^{E(T)}:
  \begin{array}{l}
    \text{the triangle inequalities are strict, and}\\
    \text{the perimeter of every face is }<2\pi
  \end{array}\right\}.
\]
The subset $\mathfrak D_+(T)\subset \mathfrak M_+(T)$ consists of spherical cone metrics for which $T$ is Delaunay.  A triangulation may have loops and multiple edges, as in \cite{IzmestievProsanovWu24}; statements about a face or the star of a vertex are understood in the abstract triangulation.

The spherical discrete conformal class of $d$ is denoted by $\mathfrak C_+(d)$.  If two metrics $d,d'\in \mathfrak D_+(T)$ differ by a vertex scaling with conformal vector $a\in \R^V$, then
\begin{equation}\label{eq:vertex-scaling-basic}
  \sin\frac{\ell_{vw}(d')}{2}
  =\exp\!\left(\frac{a_v+a_w}{2}\right)
   \sin\frac{\ell_{vw}(d)}{2}
\end{equation}
for every edge $vw\in E(T)$. Finitely many edge flips satisfying the Delaunay condition are allowed \cite[Definition 2.8]{IzmestievProsanovWu24}.

\subsection{The parameter domain}

Izmestiev-Prosanov-Wu construct an open set $U'(d)\subset\R^V$ parameterizing the spherical discrete conformal class and let $U(d)$ be the connected component containing the origin.  For a triangulation $T$, the corresponding Delaunay cell is denoted by $U(d,T)$.  
For later use, set
\[
  \Kpos:=\left\{\kappa\in(0,2\pi)^V:
  \sum_{v\in V}\kappa_v<4\pi,
  \quad 2\kappa_v<\sum_{w\in V}\kappa_w\ \text{for every }v\in V
  \right\}.
\]

We use the following results from Izmestiev-Prosanov-Wu~\cite{IzmestievProsanovWu24}.

\begin{proposition}[Izmestiev-Prosanov-Wu \cite{IzmestievProsanovWu24}]\label{prop:known-results}
Let $d\in \mathfrak M_{+}(T)$.
\begin{enumerate}
\item \cite[Lemma 2.28, Lemma 2.32]{IzmestievProsanovWu24} The map
\[
  \Psi_d:U(d)\longrightarrow \mathfrak C_+(d),\qquad u\longmapsto u*d,
\]
is a homeomorphism, and the curvature map \eqref{eq:kappa-map} is $C^1$ .

\item \cite[Corollary 2.30, Section 2.1]{IzmestievProsanovWu24} Only finitely many Delaunay cells occur:
\begin{equation}\label{eq:finite-triangulations}
  \#\{T:U(d,T)\ne\varnothing\}<\infty.
\end{equation}
Moreover, each $\mathfrak D_+(T)$ is closed in $\mathfrak M_+(T)$.

\item \cite[Remark 2.29]{IzmestievProsanovWu24} If $u,u'\in U(d,T)$ and $e=vw\in E(T)$, then
\begin{equation}\label{eq:scaling-parameter}
  \sin\frac{\ell_e(u')}{2}
  =\exp\!\left(\frac{u_v+u_w-u'_v-u'_w}{2}\right)
   \sin\frac{\ell_e(u)}{2}.
\end{equation}

\item \cite[Theorem 1.5, Lemma 3.1]{IzmestievProsanovWu24} On the positive-curvature domain $\UC(d):=\kappa^{-1}((0,2\pi)^{V})$, the map $\kappa$ is a homeomorphism onto $\Kpos$, and the Jacobian $D\kappa$ is nonsingular.
\end{enumerate}
\end{proposition}

\section{Compactness lemmas}\label{sec:compactness}

In this section, we prove two compactness lemmas needed for the proof of Theorem~\ref{thm:main}.
The arguments in this section overlap with the arguments in Izmestiev-Prosanov-Wu's work \cite{IzmestievProsanovWu24}. For completeness of the proof, we include them here.

\subsection{Finite Delaunay boundary}

\begin{lemma}[Finite Delaunay boundary]\label{lem:finite-boundary}
Fix a triangulation $T$ of $(S,V)$ and let
\[
  d_m\in \mathfrak D_+(T),
  \qquad
  \ell^{(m)}\longrightarrow\ol\ell\in(0,\pi]^{E(T)}.
\]
Assume that $\ol\ell$ does not determine a point of $\mathfrak M_+(T)$.  Then
\[
  \liminf_{m\to\infty}\Area(d_m)\ge2\pi.
\]
\end{lemma}

\begin{proof}
Each face of $d_m$ is a convex spherical triangle.  Its side lengths lie in $(0,\pi)$, satisfy the strict triangle inequalities, and have perimeter $<2\pi$.  Since every coordinate of $\ol\ell$ is positive and $\ol\ell\notin \mathfrak M_+(T)$, at least one face reaches the boundary of the space of convex spherical triangles.

We shall also use the chordal transform of Izmestiev-Prosanov-Wu, illustrated in Figure~\ref{fig:chordal trans}. Let $\mathfrak M_0(T)$ denote the space of Euclidean cone metrics on $(S,V)$ realizing $T$, and let $\mathfrak D_0(T)\subset\mathfrak M_0(T)$ be the subset for which $T$ is Delaunay.  For $d_m\in\mathfrak D_+(T)$, define
$\Phi_T(d_m)\in \mathfrak D_0(T)$
by replacing every spherical edge length $\ell$ with the Euclidean chord length
$\widehat\ell=2\sin\frac{\ell}{2}$.

\begin{figure}[h]
        \centering
        \includegraphics[width=1\linewidth]{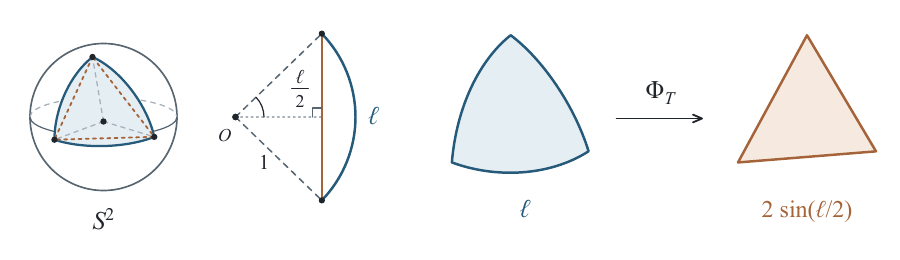}
        \caption{Chordal transform}
        \label{fig:chordal trans}
    \end{figure}
    
The image of $\Phi_T$ consists of Euclidean metrics in $\mathfrak D_0(T)$ whose face circumradii are $<1$, and the spherical Delaunay condition is equivalent to the Euclidean Delaunay condition after this transformation; see \cite[Lemma 2.11 and the discussion preceding it]{IzmestievProsanovWu24}.

There are three cases, illustrated in Figure~\ref{fig:3 cases in lem 3.1}.

\begin{figure}[h]
    \centering
    \includegraphics[width=1\linewidth]{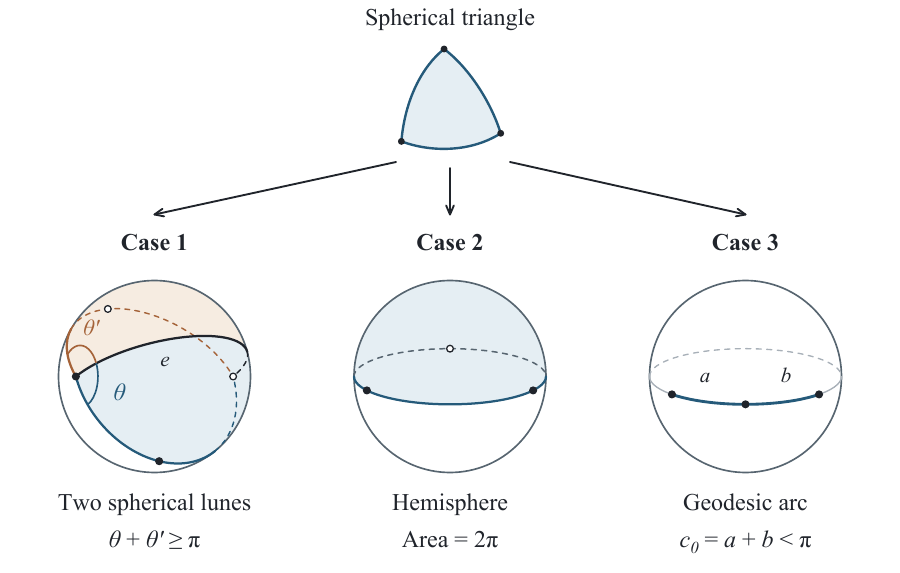}
    \caption{The three cases in Lemma 3.1}
    \label{fig:3 cases in lem 3.1}
\end{figure}

\smallskip
\noindent\textbf{Case 1: some limiting edge has length $\pi$.}

Let $e$ be such an edge, with $\ell^{(m)}(e)\to\pi$, and let $F,F'$ be the two face incidences along $e$.  If $a_m,b_m$ are the other two side lengths of $F$, then
\[
  \ell^{(m)}(e)<a_m+b_m<2\pi-\ell^{(m)}(e),
\]
which implies
\[
  a_m+b_m\longrightarrow\pi.
\]
The same argument holds for $F'$.  Thus, after developing the two faces on the unit sphere and passing to a subsequence, their limits are possibly degenerate spherical lunes with common limiting semicircle $e$.  Denote their limiting angles by $\theta,\theta'\in[0,\pi]$.


Suppose the inner angles of the spherical triangles $F$ and $F'$  are
\[
(\alpha_m,\beta_m,\gamma_m),
\qquad
(\alpha'_m,\beta'_m,\gamma'_m),
\]
where the angles $\gamma_m$ and $\gamma'_m$ are opposite to $e$.
Then
\[
(\alpha_m,\beta_m,\gamma_m)
\longrightarrow(\theta,\theta,\pi),
\qquad
(\alpha'_m,\beta'_m,\gamma'_m)
\longrightarrow(\theta',\theta',\pi).
\]
By the spherical Delaunay condition
(see Luo~\cite{Luo2006}, page 410), we have
\[
\gamma_m+\gamma'_m
\le
\alpha_m+\beta_m+\alpha'_m+\beta'_m.
\]
Passing to the limit gives
\[
2\pi\le 2\theta+2\theta',
\]
and hence $\theta+\theta'\ge\pi$. Consequently,
\[
\lim_{m\to\infty}
\bigl(\Area(F, d_m)+\Area(F', d_m)\bigr)
=2(\theta+\theta')\ge2\pi.
\]
In particular, if one face degenerates to a geodesic
arc of length $\pi$, then $\theta=0$ and
$\theta'=\pi$. Thus the adjacent face converges
to a hemisphere.

\smallskip
\noindent\textbf{Case 2: some face has limiting perimeter $2\pi$, and every limiting edge has length $<\pi$.}

Let $F$ be such a face and write its limiting side lengths as
\[
  a,b,c\in(0,\pi),
  \qquad
  a+b+c=2\pi.
\]
Let $\alpha_m,\beta_m,\gamma_m$ be the corresponding angles.  By the spherical cosine law,
\[
  \cos \alpha_m
  =\frac{\cos a_m-\cos b_m\cos c_m}{\sin b_m\sin c_m}
  \longrightarrow
  \frac{\cos(b+c)-\cos b\cos c}{\sin b\sin c}
  =-1.
\]
Thus $\alpha_m\to\pi$. Similarly, $\beta_m,\gamma_m\to\pi$.  The area formula for a spherical triangle gives
\[
  \Area(F, d_m)=\alpha_m+\beta_m+\gamma_m-\pi\longrightarrow2\pi.
\]
Hence
\[
  \liminf_{m\to\infty}\Area(d_m)\ge2\pi.
\]

\smallskip
\noindent\textbf{Case 3: a face degenerates to a geodesic arc of length $<\pi$.}

Let $F$ be such a face.  Denote its side lengths by $a_m,b_m,c_m$, where $c_m$ is the longest one, and suppose
\[
  a_m\to a,
  \qquad
  b_m\to b,
  \qquad
  c_m\to c_0=a+b<\pi.
\]
For the chordal Euclidean triangle $\widehat F=\Phi_T(F)$, put
\[
  \widehat a_m=2\sin\frac{a_m}{2},\qquad
  \widehat b_m=2\sin\frac{b_m}{2},\qquad
  \widehat c_m=2\sin\frac{c_m}{2},
\]
and let $\alpha_m$ be the angle in the Euclidean triangle $\widehat F$ opposite the edge with length $\widehat c_m$.  The Euclidean cosine law gives
\[
  \cos\alpha_m
  =\frac{\widehat a_m^2+\widehat b_m^2-\widehat c_m^2}
         {2\widehat a_m\widehat b_m}
  \longrightarrow
  -\cos\frac{c_0}{2}.
\]
Hence
\[
  \alpha_m\longrightarrow\pi-\frac{c_0}{2}.
\]

Let $F'$ be the face adjacent to $F$ along the side with length $c_m$, let $\beta_m$ be the opposite angle in $\widehat F'$, and let $R'_m<1$ be its Euclidean circumradius.  The Euclidean Delaunay condition gives
\[
  \alpha_m+\beta_m\le\pi,
\]
so
\begin{equation}\label{Eq: limit of beta}
   \lim_{m\to\infty}\beta_m\le\frac{c_0}{2}<\frac{\pi}{2}.  
\end{equation}
By the Euclidean sine law,
\[
  2\sin\frac{c_m}{2}=2R'_m\sin\beta_m.
\]
Since $c_0<\pi$, we have $\beta_m<\pi/2$ for large $m$ by (\ref{Eq: limit of beta}).  As $R'_m<1$, we have
\[
  \sin\beta_m>\sin\frac{c_m}{2},
\]
hence $\beta_m>c_m/2$. As a consequence, we have $\lim_{m\rightarrow \infty} \beta_m\geq \frac{c_0}{2}$. Therefore,
\[
  \beta_m\longrightarrow\frac{c_0}{2},
  \qquad
  R'_m=\frac{\sin(c_m/2)}{\sin\beta_m}\longrightarrow1.
\]
If $h'_m$ denotes the distance from the plane of $\widehat F'$ to the origin, then
\[
  (R'_m)^2+(h'_m)^2=1,
\]
so $h'_m\to0$.  Thus the circumcircle tends to a great circle and the spherical circumdisk tends to a hemisphere.

Consequently $F'$ reaches the boundary of the spherical triangle domain as $m\rightarrow \infty$.  If one of its limiting edges has length $\pi$, we are in Case~1; if its limiting perimeter is $2\pi$, we are in Case~2.  Otherwise, $F'$ again degenerates to a geodesic arc of length $<\pi$.

In the last case, the common side cannot be the longest side of the new arc.  Otherwise, the preceding cosine-law computation would give
\[
  \beta_m\longrightarrow\pi-\frac{c_0}{2},
\]
contrary to $\beta_m\to c_0/2$.  Hence the new arc has a longest side of limiting length $c_1>c_0$.  Repeating the argument, as long as neither Case~1 nor Case~2 occurs, we obtain
\[
  c_0<c_1<c_2<\cdots,
\]
where every $c_j$ is the limiting length of an edge of the fixed triangulation $T$.  Since $E(T)$ is finite, some edge would eventually repeat, giving
\[
  c_i<c_{i+1}<\cdots<c_j=c_i,
\]
a contradiction.  Thus the process stops after finitely many steps, necessarily in Case~1 or Case~2.  In either case
\[
  \liminf_{m\to\infty}\Area(d_m)\ge2\pi.
\]

The three cases exhaust the boundary degenerations.  This proves the lemma.
\end{proof}


The constant $2\pi$ is the area of a hemisphere. Thus Lemma \ref{lem:finite-boundary}, together with the discrete Gauss-Bonnet formula, is the source of the lower bound $\sum_v\kappa_v>2\pi$ in the definition of $K_S$.

\subsection{Escaping sequences}

\begin{lemma}[Escaping Lemma]\label{lem:escaping}
Let $u^{(m)}\in U(d)$ be a sequence leaving every compact subset of $U(d)$, and assume
\[
  \kappa\bigl(u^{(m)}\bigr)\longrightarrow\ol\kappa
  \in(-\infty,2\pi)^V.
\]
After passing to a subsequence, at least one of the following alternatives occurs:
\begin{enumerate}
\item[(E1)] $\displaystyle 2\ol\kappa_v\ge \sum_{w\in V}\ol\kappa_w$ for some $v\in V$; 
\item[(E2)] $\displaystyle \sum_{v\in V}\ol\kappa_v=4\pi$;
\item[(E3)] $\displaystyle \sum_{v\in V}\ol\kappa_v\le2\pi$.
\end{enumerate}
Thus every interior curvature limit of an escaping sequence lies in one of the three loci above.
\end{lemma}

\begin{proof}
Set $d_m:=u^{(m)}*d$.
By \eqref{eq:finite-triangulations}, after passing to a subsequence, there is a triangulation $T$ such that
\[
  d_m\in \mathfrak D_+(T)
  \qquad\text{for every }m.
\]
Because $V$ and $E(T)$ are finite, we may pass to a further subsequence so that
\[
  u_v^{(m)}\longrightarrow\ol u_v\in[-\infty,+\infty]
  \qquad (v\in V),
\]
and
\[
  \ell_e^{(m)}\longrightarrow\ol\ell_e\in[0,\pi]
  \qquad (e\in E(T)),
\]
and all face angles converge.

Taking $d_1$ as reference, the edge-length relation \eqref{eq:scaling-parameter} gives
\begin{equation}\label{eq:edge-C}
  \sin\frac{\ell_{vw}^{(m)}}2
  =C_{vw}\exp\!\left(-\frac{u_v^{(m)}+u_w^{(m)}}2\right),
\end{equation}
where
\[
  C_{vw}:=
  \sin\frac{\ell_{vw}^{(1)}}2
  \exp\!\left(\frac{u_v^{(1)}+u_w^{(1)}}2\right)>0
\]
for every edge $vw\in E(T)$.  We distinguish the same three extended-coordinate regimes as in the proof of \cite[Lemma 3.2]{IzmestievProsanovWu24}.

\medskip
\noindent\textbf{Case A: some coordinate tends to $-\infty$.}
Assume $u_v^{(m)}\longrightarrow-\infty$
for some vertex $v$.  There is no loop based at $v$ in this case. Otherwise, for such an edge, the right-hand side of \eqref{eq:edge-C} would tend to $+\infty$. Thus, every edge incident to $v$ joins $v$ to another marked vertex.  Since the left-hand side of \eqref{eq:edge-C} is at most $1$, the quantity $u_v^{(m)}+u_w^{(m)}$ is bounded from below for every incident edge $vw$.  Therefore,
\begin{equation}\label{eq:neighbors-plus-infty}
  u_w^{(m)}\longrightarrow+\infty
  \qquad\text{for every neighbor }w\text{ of }v.
\end{equation}
Let $F=\{v,w,z\}$ be a face incident to $v$.  The edge $wz$, which is opposite $v$, satisfies
\[
  \sin\frac{\ell_{wz}^{(m)}}2
  =C_{wz}\exp\!\left(-\frac{u_w^{(m)}+u_z^{(m)}}2\right)
  \longrightarrow0,
\]
so
\begin{equation}\label{eq:opposite-zero}
  \ell_{wz}^{(m)}\longrightarrow0.
\end{equation}
There are two possibilities.

\smallskip
\noindent\emph{Case A-1: no edge incident to $v$ tends to $\pi$.}

Every angle at $v$ tends to zero.  For a fixed incident face $F=\{v,w,z\}$, write
\[
  a_m:=\ell_{wz}^{(m)},\qquad
  b_m:=\ell_{vz}^{(m)},\qquad
  c_m:=\ell_{vw}^{(m)},
\]
and let $A_m$ be the angle at $v$.  By \eqref{eq:opposite-zero}, $a_m\to0$, while the triangle inequality gives
\[
  |b_m-c_m|<a_m,
\]
so $b_m$ and $c_m$ have the same limit $L\in[0,\pi)$.  If $L>0$, the spherical cosine law gives
\[
  \cos A_m
  =\frac{\cos a_m-\cos b_m\cos c_m}{\sin b_m\sin c_m}
  \longrightarrow1,
\]
hence $A_m\to0$.

Suppose now that $L=0$.  Dividing the scaling formulas for the edges $wz$ and $vz$ gives
\begin{equation}\label{eq:ratio-small}
  \frac{\sin(a_m/2)}{\sin(b_m/2)}
  =\frac{C_{wz}}{C_{vz}}
   \exp\!\left(\frac{u_v^{(m)}-u_w^{(m)}}2\right)
  \longrightarrow0
\end{equation}
by \eqref{eq:neighbors-plus-infty}.  Thus $a_m/b_m\to0$.  Since $|b_m-c_m|<a_m$, we also have $c_m/b_m\to1$.  Expanding the spherical cosine law at the origin,
\[
  \begin{aligned}
  \cos A_m
  &=\frac{\cos a_m-\cos b_m\cos c_m}{\sin b_m\sin c_m}\\
  &=\frac{b_m^2+c_m^2-a_m^2+o(b_m^2+c_m^2+a_m^2)}
          {2b_mc_m+o(b_mc_m)}
  \longrightarrow1.
  \end{aligned}
\]
Again $A_m\to0$.

Since there are finitely many faces around $v$, the cone angle satisfies $\omega_v(d_m)\longrightarrow0$.
Consequently,
\[
  \ol\kappa_v
  =2\pi-\lim_{m\to\infty}\omega_v(d_m)
  =2\pi.
\]
This contradicts the assumption $\ol\kappa\in(-\infty,2\pi)^V$.  Hence, Case~A-1 cannot occur for an interior curvature limit.

\smallskip
\noindent\emph{Case A-2: some edge incident to $v$ tends to $\pi$.}

Suppose $\ell_{vw}^{(m)}\to\pi$.  
Consider either of the two face incidences adjacent to \(vw\), say $F=\{v,w,z\}$. Then $l^{(m)}_{wz}\to 0$  by \eqref{eq:opposite-zero}, and  the triangle inequality
\[
  \bigl|\ell_{vw}^{(m)}-\ell_{vz}^{(m)}\bigr|<\ell_{wz}^{(m)}
\]
implies $\ell_{vz}^{(m)}\to\pi$. 
Iterating this argument around the cyclic link of \(v\) gives
\begin{equation}\label{eq:all-incident-pi}
  \ell_{vx}^{(m)}\longrightarrow\pi
  \qquad\text{for every edge }vx\text{ incident to }v.
\end{equation}

Let $A_{F,m}$ be the angle of an incident face $F$ at $v$, and let $\ol A_F$ be its limit.  By \eqref{eq:opposite-zero} and \eqref{eq:all-incident-pi}, each such face converges to a spherical lune whose angle at $v$ is $\ol A_F$.  A spherical lune of angle $\theta$ has area $2\theta$; hence
\[
  \Area(F,d_m)\longrightarrow2\ol A_F.
\]
Summing over the faces incident to $v$ gives
\begin{align*}
  \lim_{m\to\infty}\sum_{F\ni v}\Area(F,d_m)
  =2\sum_{F\ni v}\ol A_F
  =2\lim_{m\to\infty}\omega_v(d_m)
  =2(2\pi-\ol\kappa_v).
\end{align*}
Since all remaining face areas are nonnegative, we have
\begin{equation}\label{eq:area-A2}
  \lim_{m\to\infty}\Area(d_m)\ge 2(2\pi-\ol\kappa_v).
\end{equation}
The limit of the total area exists because of \eqref{eq:GB} and the assumed convergence of $\kappa(d_m)$.  Applying the discrete Gauss-Bonnet formula \eqref{eq:GB} and \eqref{eq:area-A2}, we have
\begin{align*}
  \sum_{w\ne v}\ol\kappa_w-\ol\kappa_v
  &=\sum_{w\in V}\ol\kappa_w-2\ol\kappa_v\\
  &=4\pi-\lim_{m\to\infty}\Area(d_m)-2\ol\kappa_v\\
  &\le0.
\end{align*}
Thus (E1) holds.

\medskip
\noindent\textbf{Case B: no coordinate tends to $-\infty$, but at least one tends to $+\infty$.}
Every $\ol u_v$ is now finite or $+\infty$.  Let
\[
  I:=\{v\in V:\ol u_v=+\infty\},
  \qquad
  F:=V\setminus I.
\]
The set $I$ is nonempty by the assumption.

If $F=\varnothing$, then both endpoints of every edge have parameter tending to $+\infty$, and \eqref{eq:edge-C} gives $\ell_e^{(m)}\to0$ for every edge.  Hence every face area tends to zero and
 $\Area(d_m)\longrightarrow0$.

Assume next that both $F$ and $I$ are nonempty.  Since the $1$-skeleton of $T$ is connected, there exists an edge that joins a vertex $w\in F$ to a vertex $v\in I$.  Consider the neighbors of $w$ in cyclic order.  If a face $\{w,v,x\}$ had $x\in F$, then the two edges incident to $v$ would tend to zero by \eqref{eq:edge-C}, while the edge $wx$ would have a positive limit because both $u_w^{(m)}$ and $u_x^{(m)}$ converge to finite numbers.  This contradicts the triangle inequality.  Therefore, $x\in I$.  Repeating the same argument around the star of $w$, every neighbor of $w$ lies in $I$.

Therefore, there is no edge with both endpoints in $F$.  Otherwise, along a shortest edge path from such an edge to $I$, the last vertex in $F$ would be adjacent both to a vertex of $F$ and to a vertex of $I$, contrary to the arguments in the preceding paragraph.  Hence, every edge of $T$ has at least one endpoint in $I$, and \eqref{eq:edge-C} again yields
\[
  \ell^{(m)}(e)\longrightarrow0
  \qquad\text{for every }e\in E(T),
\]
which implies
\[
  \Area(d_m)\longrightarrow0.
\]
By the discrete Gauss-Bonnet formula \eqref{eq:GB}, we have
\[
  \sum_{v\in V}\ol\kappa_v=4\pi,
\]
which is (E2).

\medskip
\noindent\textbf{Case C: all coordinates converge to finite real numbers.}
Assume $\ol u_v\in\R$
  for every $v\in V$.
Then \eqref{eq:edge-C} implies 
  $\ol\ell_e>0$ for every $e\in E(T)$.
Suppose that $\ol\ell\in \mathfrak M_+(T)$.  Since $\mathfrak D_+(T)$ is closed in $\mathfrak M_+(T)$, the metrics $d_m$ converge to a metric $\ol d\in \mathfrak D_+(T)$.  Passing to the limit in the scaling identity with the fixed reference metric $d_1$ gives, for every edge $e=vw$,
\[
  \sin\frac{\ell_e(\ol d)}2
  =\exp\!\left(
     \frac{u_v^{(1)}-\ol u_v+u_w^{(1)}-\ol u_w}{2}
   \right)
   \sin\frac{\ell_e(d_1)}2.
\]
Thus $\ol d$ and $d_1$ are discretely conformally equivalent with respect to the same Delaunay triangulation $T$; see \cite[Definition 2.8]{IzmestievProsanovWu24}.  Since $d_1\in \mathfrak C_+(d)$, it follows that $\ol d\in \mathfrak C_+(d)$.

Now $d_m=\Psi_d(u^{(m)})\to\ol d$ in $\mathfrak C_+(d)$.  The inverse of the homeomorphism $\Psi_d:U(d)\to \mathfrak C_+(d)$ is continuous, so
\[
  u^{(m)}=\Psi_d^{-1}(d_m)
  \longrightarrow \Psi_d^{-1}(\ol d)\in U(d).
\]
This contradicts the assumption that $u^{(m)}$ leaves every compact subset of $U(d)$.

Therefore,
$\ol\ell\in(0,\pi]^{E(T)}\setminus \mathfrak M_+(T)$,
and Lemma~\ref{lem:finite-boundary} gives
\[
  \liminf_{m\to\infty}\Area(d_m)\ge2\pi.
\]
Because $\kappa(d_m)\to\ol\kappa$, the discrete Gauss-Bonnet formula \eqref{eq:GB} implies that the areas actually converge
$\Area(d_m)\longrightarrow 4\pi-\sum_{v\in V}\ol\kappa_v$.
Therefore, $4\pi-\sum_{v\in V}\ol\kappa_v\ge2\pi$,
that is,
  $\sum_{v\in V}\ol\kappa_v\le2\pi$.
Thus (E3) holds, and the proof is complete.
\end{proof}

The lemma gives only a necessary condition for escape; no converse is asserted.  This is sufficient for the properness argument below.

\section{Proof of the main theorem}\label{sec:proof-main}

\begin{proof}[Proof of Theorem~\ref{thm:main}]
We first prove properness.  Let $L\Subset\KS$ be compact and let $u^{(m)}$ be a sequence in $\kappa^{-1}(L)$.
Suppose that no subsequence converges in $U(d)$.  Then $u^{(m)}$ leaves every compact subset of $U(d)$; otherwise some compact set would contain infinitely many terms and hence a convergent subsequence.  After passing to a subsequence, assume
\[
  \kappa\bigl(u^{(m)}\bigr)\longrightarrow\ol\kappa\in L.
\]
Set $S(\ol\kappa)=\sum_{v\in V}\ol\kappa_v$ and
apply Lemma~\ref{lem:escaping} to the sequence $\{u^{(m)}\}$.  
Alternative (E1) contradicts the inequalities $2\ol\kappa_v<S(\ol\kappa)$ defining $K_S$.  Alternative (E2) contradicts $S(\ol\kappa)<4\pi$, and (E3) contradicts $S(\ol\kappa)>2\pi$.  Thus every sequence in $\kappa^{-1}(L)$ has a convergent subsequence in $U(d)$.  Since $\kappa^{-1}(L)$ is closed in $U(d)$, it is compact.  Hence, the curvature map
\[
  \kappa:\XS(d)\longrightarrow\KS
\]
is proper.

It remains to compute the degree.  Set $\kappa^*:=\left(\frac{3\pi}{n},\ldots,\frac{3\pi}{n}\right)$.
Then
\[
  \sum_v\kappa_v^*=3\pi,
  \qquad
  2\kappa_v^*=\frac{6\pi}{n}<3\pi
\]
for every $v\in V$, so $\kappa^*\in\KS\cap\Kpos$.  By Proposition~\ref{prop:known-results}(4), $\kappa^*$ has a unique preimage $u^*\in\UC(d)$ and $D\kappa(u^*)$ is nonsingular.

Fix an ordering of $V$ and use the induced orientation of $\R^V$.  Properness gives a degree for the restricted map.  Since $K_S$ is convex, this degree is computed at the regular value $\kappa^*$:
\[
  \deg\bigl(\kappa|_{\XS(d)}\bigr)
  =\operatorname{sgn}\det D\kappa(u^*)
  =:\eps_d\in\{+1,-1\}.
\]
In particular, the degree is nonzero.  The map is therefore onto.

Finally, let $\ol\kappa\in K_S$ be regular. $\kappa^{-1}(\ol\kappa)$ is discrete, and properness makes it compact; hence it is finite.  The degree formula gives
\[
  \sum_{u\in\kappa^{-1}(\ol\kappa)}
  \operatorname{sgn}\det D\kappa(u)=\eps_d.
\]
Each summand is $\pm1$, while the sum is $\pm1$.  The number of terms is therefore odd.  This proves all assertions.
\end{proof}

\end{document}